\documentclass[12pt]{amsart}
\usepackage{amsfonts}
\usepackage{amssymb}
\input epsf
\usepackage{amsmath}
\usepackage{eucal}
\usepackage{graphicx}
\usepackage{amscd}

\usepackage{amscd}
\begin{document}

\title{\large  Matricial formulae for partitions }

\author[F. Aicardi]{F. Aicardi}

\address[F.~Aicardi]{SISTIANA 56 PR (Trieste),  Italy}
\email{aicardi@sissa.it }
%
%
%

%
\begin{abstract}
  The  exponential
of  the  triangular matrix  whose entries in the diagonal  at  distance $n$  from the principal diagonal
are all equal  to the sum of the inverse of the  divisors  of $n$ is  the triangular  matrix whose  entries in the diagonal at   distance $n$  from the
principal  diagonal  are all equal to the  number  of partitions of $n$.  A similar result is true for any
pair  of sequences satisfying  a special recurrence.

\end{abstract}
\maketitle

Let  $\Delta$ be  the infinite
triangular matrix having zeroes on   and below the
main diagonal,  and whose values on each parallel to the main diagonal  are all equal to:
\[ {\Delta}_{i, i+n}= \overline\sigma(n)  \ \ \ \ \  (i \ge 0,    \ \ \  n>0),  \]
where $\overline\sigma(n)$ is the sum of the inverses of the divisors of $n$, for $n>0$.

{\bf  Proposition 1.} {\it The
matrix
\[  {\bf P}=\exp (\Delta) \]
is triangular,  with zero below the diagonal,  $1$ on the diagonal,
and  its elements on the  parallel to the diagonal
at distance $n$ from the diagonal are  equal
to  the number $p(n)$ of partitions of the integer $n$
  (being $p(0)=1$): }
\[ {\bf P}_{i, i+n}= p(n)  \ \ \ \ \  (i\ge 0, \ \  n\ge 0).\]

\[ {\Delta}=\left(  \begin{array}{cccccccc}
     0 & \frac{1}{1} & \frac {3}{2} & \frac{4}{3} &  \frac{7}{4} & \frac{6}{5} & \frac{12}{6}& \dots  \\[5 pt]
     0 & 0 & \frac{1}{1} & \frac{3}{2} & \frac{4}{3} & \frac{7}{4} &\frac{6}{5} & \dots  \\[5 pt]
    0 & 0  & 0 & \frac{1}{1}& \frac{3}{2} & \frac{4}{3} & \frac{7}{4} & \dots \\[5 pt]
    0& 0& 0 & 0  & \frac{1}{1} &  \frac{3}{2} & \frac{4}{3}& \dots \\[5 pt]
    0& 0& 0& 0 & 0 & \frac{1}{1} &\frac{3}{2} & \dots \\[5 pt]
     0& 0& 0& 0& 0 & 0 & \frac{1}{1}&  \dots \\[5 pt]
    0& 0& 0 & 0& 0& 0 &0 & \dots  \\[5 pt]
   \vdots & \vdots  & \vdots & \vdots &  \vdots& \vdots & \vdots & \ddots
\end{array}   \right) \ \ \ \ \
 {\bf P}=\left(  \begin{array}{cccccccc}
     1 & 1 & 2 & 3 & 5 & 7 &11 & \dots \\[5 pt]
     0 & 1 & 1 & 2 & 3 & 5 &5 & \dots \\[5 pt]
        0 & 0& 1 & 1 & 2 & 3 &5 & \dots  \\[5 pt]
      0 & 0 & 0 & 1 & 1 & 2 & 3& \dots  \\[5 pt]
      0 & 0 & 0 & 0 & 1 & 1 & 2 & \dots  \\[5 pt]
      0 & 0 & 0 & 0 & 0 & 1 & 1 & \dots  \\[5 pt]
     0 & 0 & 0 & 0 & 0 & 0 &1  & \dots  \\[5 pt]
     \vdots & \vdots  & \vdots & \vdots &  \vdots& \vdots & \vdots & \ddots
     \end{array}   \right)\]

We prove in the next section a more general result.
As  byproduct,  we find a  formula  (see  eq.(\ref{eq3})) for the  number of  partitions of $n$
as sum   of the products  of the function  $\overline \sigma$ over all  elements of
each partition, and  a formula for $\overline \sigma$  as combination     of the products of the
number  of partitions  over all elements of each partition (see  eq.(\ref{eq7})).

{\bf Generalisation of Proposition 1.}

Let $\{s_n\}_{n=0}^\infty $ be any infinite sequence,  satisfying $s_0=0$.

Let $t_n$ be the infinite sequence obtained by the following recurrence:
\begin{equation}\label{rec} t_0=1,  \quad \quad  t_n=\frac{1}{n}\sum_{k=0}^{n-1} s_{n-k}t_k  \end{equation}

 Consider  the infinite
triangular  matrix ${\bf S}$   having zeroes on the diagonal and below the
diagonal,  and whose value on each parallel to the diagonal  and at distance $n$ from the diagonal  are all equal to:
\[ {\bf S}_{i, i+n}= s_{n}/n  \ \ \ \ \  (i\ge 0,   \ \ \  n > 0 ). \]

{\bf Theorem 1.} {\it
The matrix
\[  {\bf T}=\exp ({\bf S}) \]
is triangular,  with zero below the diagonal,  $1$ on the diagonal,
and the values of its elements on each parallel to the diagonal
 at distance $n$ from the diagonal are equal
to   $t_n$: }
\[ {\bf T}_{i, i+n}= t_{n}  \ \ \ \ \  (i\ge 0,  \ \  n \ge 0 ).\]

{\bf  Example 1.} Proposition 1 is  a particular  case of Theorem 1.  Indeed,  $\overline\sigma(k)=\sigma(k)/k$, $\sigma(k)$ being  the sum of all divisors  of $k$.   The recurrence
\[  np(n)=\sum_{k=0}^{n-1} \sigma(n-k)p(k)\equiv \sum_{h=1}^n \sigma(h)p(n-h)  \]
is proved, for instance in [1].  (A proof can be done  in a similar way as we do  for the recurrence of  the next example).

{\bf  Example 2.} Let $p^o(n)$  be  the  number of  partitions  of  $n$  into {\bf odd}
integers,

\[ \begin{array}{c||c|c|c|c|c|c|c|c|c|c|c|c|c|c|c|c|c|c|c|c|c}
   n & 1 & 2 & 3 & 4 & 5 & 6 & 7  & 8 & 9 & 10  & 11& 12 & 13 & 14 & 15  & 16 &17 &18 &19 & 20& \dots \\  \hline
    p^o(n)  &  1& 1& 2& 2& 3& 4& 5& 6& 8& 10& 12 & 15 & 18 & 22 & 27& 32& 38& 46& 54& 64& \dots  \end{array}\]
and    $\sigma^o_k$  be the sum of {\bf odd} divisors  of $k$:

\[ \begin{array}{c||c|c|c|c|c|c|c|c|c|c|c|c|c|c|c|c|c|c|c|c|c}
   k & 1 & 2 & 3 & 4 & 5 & 6 & 7  & 8 & 9 & 10  & 11& 12 & 13 & 14 & 15& 16 &17 &18 &19 & 20 & \dots \\ [2 pt ] \hline
    \sigma^o(k)   &  1&  1&  4&  1&  6&  4&  8&  1&  13&  6 &  12&  4&  14&  8&  24&  1&  18&  13&  20&  6 & \dots  \\  \end{array}\]

{\bf Lemma. } \[   np^o(n)=\sum_{k=0}^{n-1} \sigma^o(n-k)p^o(k)\equiv \sum_{h=1}^n \sigma^o(h)p^o(n-h).   \]

{\it Proof of Lemma.}  Consider all  partitions of $n$  into odd  integers
\[ n= o_1+o_2+ \dots+ o_r , \quad \quad   o_i \quad  {\rm odd}, \]
and   sum all these  equalities:  we obtain at left $np^o(n)$, and at right
\[   \sum_{i=1}^n  o_i \mu_i,  \quad \quad  \mu_i\ge 0,   \]
$\mu_i$  being  the number of  times the odd $o_i$ appears in total  inside all partitions of  $n$
into  odd  integers.   We have
\[ \mu_i= \sum_{m=1}^{n/o_i} p^o(n-m o_i).  \]
Indeed,  there are possibly several  copies of $o_i$ in a partition:  $p^o(n-o_i)$ is the number of those
containing at least one copy of $o_i$, and then we count in this summand  the first  copy of the $o_i$
in all partitions; $p^o(n-2o_i)$
is the number of partitions containing  at least 2 copies of $o_i$, and then we count in this summand the
second copy of  $o_i$ in all partitions, and so on.
We obtain:
\[  n p^o(n)=  \sum_{i=1}^n  o_i \sum_{m=1}^{n/o_i} p^o(n-m o_i). \]
Writing  $m o_i=h$,  we  get
\[  n p^o(n)=  \sum_{h=1}^n   \sum_{o_i | h} o_i    p^o(n-h), \]
the second sum being extended to all  odd  integers dividing  $h$.  Then
\[  p^o(n)= \frac{1}{n} \sum_{h=1}^n   \sigma^o(h) p^o(n-h) \equiv \frac{1}{n} \sum_{k=0}^{n-1}   \sigma^o(n-k) p^o(k).  \]
\hfill $\square$

 Define the  matrix $\Sigma^o$  whose elements are zero on and below the diagonal and
 \[  \Sigma^o_{i, i+n}= \sigma^o(n)/n   \quad \quad (i\ge 0, \quad  n>0). \]
By Theorem 1,   the matrix ${\bf P}^o=\exp(\Sigma^o)$ is triangular  and  satisfies:
\[   {\bf P}^o_{i, i+n}= p^o(n),    \quad \quad (i\ge 0, \quad n\ge 0 ).  \]

{\bf  Example 3.}  Let $p^{\not q}(n)$  be  the  number of  partitions  of  $n$  into
integers  non  divisible  by  the prime $q$,
  and    $\sigma^{\not q}_k$  be the sum of  divisors  of $k$  non divisible by $q$:

Like in the  above lemma, we prove

 \[   np^{\not q}(n)=  \sum_{h=1}^n \sigma^{\not q}(h)p^{\not q}(n-h).   \]

 Define the  matrix $\Sigma^{\not q}$  whose elements are zero on and below the diagonal and
 \[  \Sigma^{\not q}_{i, i+n}= \sigma^{\not q}(n)/n   \quad \quad (i\ge 0, \quad  n>0). \]
By Theorem 1,   the matrix ${\bf P}^{\not q}=\exp(\Sigma^{\not q})$ is triangular  and  satisfies:
\[   {\bf P}^{\not q}_{i, i+n}= p^{\not q }(n)     \quad \quad (i\ge 0, \quad  n\ge 0 ).  \]

{\bf Proof of Theorem 1}

 The exponential of ${\bf S}$ is given by:  \begin{equation} \label{eq1} \exp({\bf S})=
{\bf E} + {\bf S}+\frac{1}{2!}{\bf S}^2+ \frac{1}{3!}{\bf
S}^3+ \frac{1}{4!}{\bf S}^4 +  \dots   \end{equation} where
$\bf E$ is the identity matrix.

Consider the sequence $\tau_n$,  defined as the   first   row  of the  matrix ${\bf T}=exp(\bf S)$.
We  will  prove  that  ${\bf T}_{i,i+n}=\tau_{n}$,  and  $\tau_n$ satisfies the recurrence  (\ref{rec}).

We  use in the following the notation  $s(n)$  instead of $s_n$  in order to avoid many subscript indices.  We thus have:
\[ {\bf S}_{i,k}=s(k-i)/(k-i), \quad \quad {\bf S}_{0,k}=s(k)/k.  \]

We calculate the powers of ${\bf S}$.

 The square  ${\bf S}^2$  is equal to:
\[   {\bf S}^2_{i,j}=  \sum_{k=0}^n {\bf S}_{i,k} {\bf S}_{k,j}= \sum_{i<k<j} \frac{s(k-i)s(j-k)}{(k-i)(j-k)},   \]
being $S_{i,k}=0$ for $k\le i $.  Hence ${\bf S}^2_{i,j}=0$ for $j<i+2$.

The matrix ${\bf S}^3$  is equal to:
\[ {\bf S}^3_{i,j}=  \sum_{k=0}^n {\bf S}_{i,k} {\bf S}^2_{k,j}=  \sum_{k=0}^n {\bf S}_{i,k}\sum_{h=0}^n {\bf S}_{k,h}{\bf S}_{h,j}= \sum_{i<k<h<j} \frac{s(k-i)s(h-k)s(j-h)}{(k-i)(h-k)(j-h)},   \]
being ${\bf S}^2_{i,j}=0$ for $j\le i+1$.   Hence ${\bf S}^3_{i,j}=0$ for $j< i+3$.

We  obtain, for the $r$-th power of ${\bf S}$:
 \[ {\bf S}^r_{i,j}=  \sum_{i=k_0<k_1<k_2< \dots < k_r=j} \ \  \prod_{m=1}^r \frac{s(k_m-k_{m-1})}{(k_m-k_{m-1})},    \]
which is zero if $j < i+r$.

The first  elements of the sequence  $\tau_0,\dots,\tau_n$ are   the first elements of the first   row  of the  matrix $\sum_{r=0}^n \frac{1}{r!}({\bf S})^r$,  since    the elements  at distance lower  than $m$  over  the diagonal of  $ {\bf S}^{m}$ are  zero, as we have   seen.

The   element ${\bf T}_{0,n}$  of the first  row  of the  matrix  ${\bf E}+\sum_{r=1}^n \frac{1}{r!}{\bf S}^r$ are therefore, for $n>0$
\[ \tau_n=  \left\{{\bf E}+\sum_{r=1}^n \frac{1}{r!} {\bf S}^r \right\}_{0,n}=  \sum_{r=1}^n \frac{1}{r!} \sum_{0<k_1<k_2< \dots < k_r=n} \ \ \prod_{m=1}^r \frac {s(k_m-k_{m-1})}{(k_m-k_{m-1})},     \]
and $\tau_0=1$.
 For every $r=1\dots n$, the $r$ positive integers $k_m$  satisfy
 \[ \sum_{m=1}^n (k_m-k_{m-1})=n; \]  hence
 the  sum is extended  to  all  partitions of $n$ containing exactly $r$  elements $n=h_1,\dots, h_r$. However,
  for every  partition, the same  product   $\prod_{m=1}^r s(h_m)/(h_m)$  corresponds to a number of cases
  with different $k_m$ depending on the ordering of the  $h_m$.   This number of different cases is equal to
  \[  r!/\prod_{m=1}^r \rho_m !, \]  where $\rho_m$  is the multiplicity of $h_m$ in the partition.

Hence we obtain

{\bf  Proposition 2.} { \it The  elements $\tau_n$, for $n>0$,  are equal to:

\begin{equation}\label{eq3}  \tau_n=\sum_{r=1}^n  \sum_{X_r(n)} \prod_{m=1}^r \frac{s(h_m)}{h_m\rho_m !}  \end{equation}

where  $X_r(n)$  is  the set  of partitions of $n$ into exactly $r$ summands $h_m$,
and $\rho_m$  is  the multiplicity of  $h_m$. }

The formula expressing the elements of  ${\bf S}^r$ shows that ${\bf S}^r_{i,j}$  depends only on the difference
$(j-i)$,  and therefore the elements  ${\bf T}_{i,j}$ are equal on each  parallel at distance $r$  from the   diagonal, where  $j=i+r$.

The proof of  Theorem 2  is now completed proving  that  formula  (\ref{eq3}) satisfies the recurrence:
\begin{equation}\label{eq4}  \tau_N=  \frac{1}{N}\sum_{n=0}^{N-1} s(N-n)\tau_{n}.\end{equation}

We write the left member  of (\ref{eq4})  by substituting  the terms $\tau_{n}$  given  by formula (\ref{eq3}):
\begin{equation}\label{eq5} \frac{1}{N} \left( s(N) + \sum_{n=1}^{N-1} s(N-n) \sum_{r=1}^n  \sum_{X_r(n)} \prod_{m=1}^r \frac{s(h_m)}{h_m\rho_m !}\right). \end{equation}

Each summand in the  expression  from any  $\tau_n$, $0<n<N$,  is the product of  $r$  factors $s(h_m)$, such that
$\sum_{m=1}^r h_m \rho(m)=n$. In  expression (\ref{eq5}) it is multiplied by $s(N-n)$ and becomes therefore the  product  of $r+1$ factors $s(h_m)$  such that their sum is $N$,  as in the expression of $\tau_N$  by  (\ref{eq3}).

Consider now any one  of such terms, containing  $\prod_{m=1}^{r+1} s(h_n)$, in which  one  factor,  say  $s(h_{m^*})$, is $s(N-n)$.

If $h_{m^*}$  is different from the other $r$ $h_m$, the  contribution  of   the corresponding term is equal to
\[ \frac{1}{N}\frac{1}{\prod_{m=1}^{r} h_m \rho_m!}= \frac{1}{N}\frac{1}{\prod_{m=1}^{r+1} h_m \rho_m! /h_{m^*}}. \]

If $h_{m^*}$  is equal to one of the other $h_m$, i.e.,  in the product  has multiplicity $\rho_m>1$, the  coefficient of the corresponding term is equal to
\[\frac{1}{N}\frac{1}  {\prod_{m=1}^{r+1} h_m \rho_m!/(h_{m^*}\rho_{m^*})}, \]
since $s(h_{m^*=N-n})$ had multiplicity  $\rho_{m^*}-1$ in  the original  term  in  $\tau_n$.

If  among the  $r+1$  integers $h_m$, $d\le r$  integers $h_m^*$ are  different,  in correspondence of the same product we are considering,  there will be in  $t_N$ the  contribution  of $d$  terms, namely  for every  value of  $n$  such that  $N-n=h_{m^*}$.

  The coefficient  of the terms
$ \prod_{m=1}^{r+1} s(h_m)$ is thus
\[  \frac{1}{N}\sum_{k=1}^d  \frac{1}{\prod_{m=1}^{r+1} h_m \rho_m!/(h_{m^*_k} \rho_{m^*_k})},  \]
which is equal to
\[ \frac{ \sum_{k=1}^d  h_{m^*_k} \rho_{m^*_k}}{N \prod_{m=2}^{r+1}h_m \rho_m!}= \frac{N}{N\prod_{m=1}^{r+1} h_m \rho_m!}=\frac{1}{\prod_{m=1}^{r+1}h_m \rho_m!}. \]

Hence we  obtain from eq. (\ref{eq5}),  denoting $q=r+1$:
\[ \tau_N= \frac{s(N)}{N}+ \sum_{r=1}^{N-1} \sum _{X_{r+1}(N)}\prod_{m=1}^{r+1}\frac{s(h_m)}{h_m \rho_m!}=\frac{s(N)}{N}+\sum_{q=2}^{N} \sum _{X_{q}(N)}\prod_{m=1}^{q} \frac {s(h_m)}{h_m  \rho_m!},\]
i.e., the expression of $\tau_N$ by eq. (\ref{eq3}).  \hfill $\square$

{\bf Inverse formula}

We  see now  how the sequence  $s_n$ can be written in terms of the sequence  $t_n$ satisfying  the recurrence (\ref{eq1}).

{\bf  Proposition 3.} {\it The following identity holds:}
\begin{equation}\label{eq7}
\frac{s_n}{n} =   \sum_{r=1}^n  (-1)^{r-1} (r-1)! \sum_{X_r(n)} \prod_{m=1}^r \frac{t_{h_m}}{\rho_m !}
\end{equation}

{\it Proof.}

From  the  identity ${\bf T}=\exp({\bf S})$ we obtain
\[    {\bf S}=\ln {\bf T}= \ln({\bf E}+ ({\bf T}-{\bf E}))=  \sum_{r\ge 1} (-1)^{r-1}\frac{1}{r}({\bf T}-{\bf E})^r.  \]

The first row  of matrix ${\bf T}-{\bf E}$   contains, but the initial  elements which is zero, the elements
of the sequence $t_n$.  Hence we calculate  the  elements  of  $({\bf T-{\bf E}})^r$,  as in the proof of
 Theorem 1, and  we  add them  with coefficients  $(-1)^{r+1}/r$, to obtain the sequence $s_n$.
 \hfill $\square$


\vskip 1 cm
\section*{References}

[1] E. Grosswald,   {\it  Topics from the theory of numbers},
Macmillan,  New York,  1966,  p. 226.

\end{document}